\documentclass[a4paper,leqno,11pt]{amsart}
\usepackage{amsmath, amssymb, amsfonts, latexsym, mdwlist, amsthm}
\usepackage{graphicx}
\usepackage{enumerate}

\usepackage[all]{xy}

\usepackage[colorlinks=true, citecolor=blue, urlcolor=blue, linkcolor=blue%, pagebackref
]{hyperref}

\newcommand\bbA{\mathbb A}
\newcommand\CC{\mathbb C}

\newcommand\Hom{\operatorname{Hom}\nolimits}

\newcommand\im{\operatorname{im}\nolimits}

\newcommand\OOO{{\mathcal O}}
\newcommand\PP{\mathbb{P}}

\newcommand\rk{\text{rk}}

\newcommand{\actson}{\raisebox{1pt}{\scalebox{0.6}{\begin{picture}(15,5)\thinlines
\put(3,4){\oval(18,10)[r]}
\put(2,9){\vector(-1,0){0.5}}\end{picture}}
}}
\newcommand{\equi}{\: \Longleftrightarrow\: }
\renewcommand{\bar}[1]{\overline{#1}}

\newtheorem{theorem}{Theorem}

\newtheorem{prop}[theorem]{Proposition}

\newtheorem{lemma}[theorem]{Lemma}

\newtheorem*{mthm}{Main Theorem}

\theoremstyle{definition}
\newtheorem{definition}[theorem]{Definition}

\newtheorem{proposition-definition}[theorem]{Proposition-Definition}

\newtheorem{remark}[theorem]{Remark}

\nonstopmode

\begin{document}

\vspace{1cm}

\title{Unirationality of instanton moduli space for small charges}

\author{Dimitri Markushevich}
\address{\scriptsize D. Markushevich:
Univ. Lille, CNRS, UMR 8524 -- Laboratoire Paul Painlev\'e, F-59000 Lille, France}
\email{dimitri.markouchevitch@univ-lille.fr}

\author{Alexander Tikhomirov}
\address{Faculty of Mathematics\\
	National Research University  
	Higher School of Economics\\
	6 Usacheva Street\\ 
	119048 Moscow, Russia}
\email{astikhomirov@mail.ru}

\begin{abstract}
The unirationality of the moduli space of mathematical instantons on the projective 3-space is proved for charges less than or equal to 7.
%AMS Mathematics Subject Classification 14D20.
\end{abstract}

\maketitle

\thispagestyle{empty}

%\begin{abstract}
%
%2000 Mathematics Subject Classification, Primary 14M15, (Secondary 14J60, 32L05).
%\end{abstract}

A {\it mathematical instanton of charge $n\geq 1$} is a rank-2 algebraic
vector bundle $E$ on the 3-dimensional projective space $\mathbb{P}^3$ with Chern classes
$
c_1(E)=0,\ c_2(E)=n,
$
satisfying the vanishing conditions
$
h^0(E)=h^1(E(-2))=0.
$
The moduli space of these objects, denoted $I_n$, is a smooth and irreducible quasi-projective variety of dimension $8n-3$ \cite{JV,Tikhomirov2}. The problem of determining the birational type of $I_n$ is quite resistant. The rationality of $I_n$ is known for charges $n=1,2,3$ and 5 (see \cite{B1,H,ES,Ka}). In the present note we prove:

\begin{mthm}
The variety $I_n$ is unirational for $n=4, 6, 7$.
\end{mthm}

By \cite{BH}, each instanton $E$ is the middle cohomology of the following 3-term complex of sheaves, called monad:
\begin{equation}\label{monad}
0\to H_n\otimes\mathcal{O}_{\mathbb{P}^3}(-1)\overset{\alpha}\to W_{2n+2}\otimes\mathcal{O}_{\mathbb{P}^3}
\overset{\alpha^\vee}\to
H_n^\vee\otimes\mathcal{O}_{\mathbb{P}^3}(1)\to0.
\end{equation}
Here $H_n$, $W=W_{2n+2}$ are vector spaces of respective dimensions $n,2n+2$, and $W$ is endowed with a symplectic form $q$ making the monad self-dual, that is the surjection ${\alpha^\vee}$ is symplectic conjugate to the injection $\alpha$. This leads to the representation of $I_n$  as the quotient U/G, where
$U$ is the set of linear algebra data parametrizing the above monads and $G=GL(H_n)\times Sp(W,q)/\{\pm(1,1)\}$
is the group acting on $U$ in a natural way. The linear algebra data in question arise in the following way: consider the map $H^0({\alpha^\vee}):W\to H_n^\vee\otimes V^\vee$, where we introduced a 4-dimensional vector space $V$ such that $\PP^3=\PP(V)$, and denote by $\gamma=\gamma(\alpha):H_n\otimes V\to W$ the conjugate
of $H^0({\alpha^\vee})$, using $q$ to identify $W$ with $W^\vee$. Thus $H^0({\alpha^\vee})=\gamma^\vee$, and we set $A=A(\gamma):= \gamma^\vee\circ \gamma$; it is a skew-symmetric map $H_n\otimes V\to H_n^\vee\otimes V^\vee$, or equivalently an element $A\in\wedge^2(H_n^\vee\otimes V^\vee)$ which will be denoted by the same symbol. The condition that the sheafifications $\alpha=\alpha(\gamma),\ \alpha^\vee$ of $\gamma,\gamma^\vee$ form a complex, that is
$\alpha^\vee\circ\alpha=0$ in \eqref{monad}, is equivalent to $A\in \wedge^2H_n^\vee\otimes S^2V^\vee$.
\begin{theorem}[Barth--Hulek \cite{BH}] There is a $G$-invariant Zariski open set 
$$U\subset \{\gamma\in\Hom(H_n\otimes V,W)\ | \ A(\gamma)\in \wedge^2H_n^\vee\otimes S^2V^\vee\}$$
such that $I_n$ is naturally isomorphic to the quotient variety $U/G$. Moreover the quotient map $U\to I_n$ is a principal $G$-bundle.
\end{theorem}

In, \cite{B2}, Barth represented an open subset of $I_n$ as a quotient of a slice of $U$ by a smaller group.

\begin{definition}
Let $X$ be an irreducible algebraic variety, $Y\subset X$ a closed irreducible subvariety,
$G$ an algebraic group acting on $X$, and $H$ a closed subgroup of $G$.
We say that $Y$ is a $(G,H)$-slice of the action $G:X\actson$ if the following two conditions
are verified:

(i) $\overline{G\cdot Y}=X$, and

(ii) there exists an open subset $Y_0\subset Y$ such that for $y\in Y_0$, we have
$$
g\in G,\ gy\in Y\ \ \equi \ \ g\in H.
$$
\end{definition}

\begin{prop}
If $Y$ is a $(G,H)$-slice of the action $G:X\actson$, then
${\mathbf{\CC}}(X)^G\simeq {\mathbf{\CC}}(Y)^H$.
\end{prop}

To describe Barth's $(G,H)$-slice $Y$, we fix some bases in $H_n$, $V$, $W$, in such a way that
$q=\left(\begin{smallmatrix}0&1_n\\-1_n&0\end{smallmatrix}\right)\oplus
\left(\begin{smallmatrix}0&1\\-1&0\end{smallmatrix}\right)$. This allows us to represent $\gamma:H_n\otimes V\to W$ by a $4n\times(2n+2)$-matrix, which we write as a $4\times 4$ block matrix. We set 
\begin{equation}\label{Barth-matrix}
\Sigma=\left\{\left.\widetilde{\gamma}=
\left(\begin{array}{cccc}
\mathbf{0} & \mathbf{1}_n & A_1 & A_2\\
-\mathbf{1}_n & \mathbf{0} & B_1 & B_2\\
\mathbf{0} & \mathbf{0}& a^T_1 & a^T_2 \\
 \mathbf{0}&\mathbf{0} & b^T_1 & b^T_2
  \end{array}\right)\ \ \right|\ \begin{minipage}{8.5truecm}$$
  A_i,B_i\in S^2H_n^\vee,\ a_i,b_i\in H_n^\vee,$$
  $$[A_1,B_1]+a_1\wedge b_1=\mathbf{0},\ 
  [A_2,B_2]+a_2\wedge b_2=\mathbf{0},$$
  $$[A_1,B_2]+[A_2,B_1]+a_1\wedge b_2+a_2\wedge b_1=\mathbf{0}.$$
  \end{minipage}\right\}
\end{equation}
where $a_i,b_i\in H_n^\vee$ are viewed as columns, so that the transposes $a^T_i,b^T_i$ are rows of length $n$, and $a\wedge b:=ab^T-ba^T$ for two columns $a,b$ of length $n$. The locus $\Sigma$ is $H$-invariant for $H=O(n)\times SL(2)/\{\pm (1,1)\}\subset G$ with the action of
$\left(g,\left(\begin{smallmatrix}s & t \\
u & v\end{smallmatrix}\right)\right)\in O(n)\times SL(2)$
 given by
\begin{multline}\label{On-action}
(A_1,A_2,B_1,B_2,a_1,b_1,a_2,b_2)\mapsto
(gA_1g^T,gA_2g^T,gB_1g^T,gB_2g^T,\\
sga_1+ugb_1,tga_1+vgb_1,sga_2+ugb_2,tga_2+vgb_2).
\end{multline}

\begin{theorem}[Barth \cite{B2}]
The intersection $Y=\Sigma\cap U$ is a $(G,H)$-slice of the action of $G$ on $U$. The quotient map
$Y\to Y/H$ is a principal $H$-bundle, and $Y/H$ is identified with a Zariski open subset $I_n^0$ of $I_n$.
\end{theorem}

{\em The Main Theorem} is a consequence of the following lemmas:

\begin{lemma}
In the neighborhood of any point of $Y$, the equations in \eqref{Barth-matrix} in the $2n(n+3)$ matrix elements of $A_i, B_i, a_i,b_i$, represent $Y$ as a smooth transversal intersection of $\frac32n(n-1)$ hypersurfaces in the affine space $\bbA^{2n(n+3)}$.
\end{lemma}

\begin{proof}
Indeed, by Barth's theorem, $Y\to I_n^0$ is a principal $H$-bundle; as $\dim H=\frac12n(n-1)+3$ and $\dim I_n^0=8n-3$, the count of parameters shows that $Y$ is a complete intersection of the equations in  \eqref{Barth-matrix}. Moreover, by \cite{JV}, $I_n^0$, and hence $Y$ is smooth, hence the Jacobian matrix of these equations is of maximal rank, equal to the number of equations.
\end{proof}

\begin{lemma}
Let $\bar Y\subset \Sigma$ be the irreducible component of $\Sigma$ containing $Y$. Consider the linear projection
$$
\pi:(S^2H_n^\vee)^4\oplus (H_n^\vee)^4\to (S^2H_n^\vee)^2\oplus (H_n^\vee)^2, \ 
(A_i,B_i,a_i,b_i)_{i=1,2}\mapsto (A_1,A_2,a_1,a_2).
$$
Then for $4\leq n\leq 7$, the restriction of $\pi$ to $\bar Y$ is an $O(n)$-equivariant vector bundle over a nonempty Zariski open subset of $(S^2H_n^\vee)^2\oplus (H_n^\vee)^2$.
\end{lemma}

\begin{proof}
The relations on $\tilde\gamma$ in \eqref{Barth-matrix} become linear equations in $B_1,B_2,b_1,b_2$ as soon as we fix $(A_1,A_2,a_1,a_2)$. The number of equations is $\frac32n(n-1)$ and the number of variables is $n(n+3)$, so the expected dimension of the space of solutions $(B_i,b_i)$ is $n(n+3)-\frac32n(n-1)=\frac12n(9-n)$. We assume that $4\leq n\leq 7$, then the expected dimension is strictly positive, and even $\geq 7$.  Hence there is an irreducible component $\Sigma_0$ of $\Sigma$ dominating $(S^2H_n^\vee)^2\oplus (H_n^\vee)^2$.
It suffices to verify that there exists a point $y\in Y$ in which the dimension of the space of solutions is the expected one, since this implies that $\bar Y=\Sigma_0$ by the previous lemma.

We now present an explicit open condition on $\tilde \gamma$ which one may add  to the equations in \eqref{Barth-matrix} in order to guarantee that the associated monad defines an instanton:
%\begin{enumerate}[(i)]
%\item $\rk \:(\, a_1+ta_2 \ \ b_1+tb_2\,)=2$ for all $t\in \CC$;
%$\rk\left(\begin{smallmatrix}a_1&b_1\\a_2&b_2\end{smallmatrix}\right)=2$;
the sheaf map $\alpha(\tilde\gamma)^\vee:W\otimes\OOO_{\PP^3}\to H_n^\vee\otimes\OOO_{\PP^3}(1)$ is surjective. This automatically implies that the conjugate injection $\alpha(\tilde\gamma)$ is a vector bundle map,
and then the cohomology of the monad $E(\tilde\gamma):= \ker\alpha(\tilde\gamma)^\vee/\im\alpha(\tilde\gamma)$ is a rank 2 vector bundle; the vanishing conditions from the definition of an instanton and the values of Chern classes easily follow from the monad. One can also easily verify that, assuming $(A_i, a_i)$ generic, the latter open condition is a consequence of the following one:
$$
\rk \:(\, a_1+ta_2 \ \ b_1+tb_2\,)=2 \mbox{\ for all\ } t\in \CC.
$$
The verification of this condition is a problem in linear algebra. Using algebra computing system Macaulay2
\cite {M2}, we checked that a random explicit solution $B_1,B_2,b_1,b_2$ for a randomly chosen quadruple $(A_1,A_2,a_1,a_2)$ satisfies the latter condition, thus $\bar Y=\Sigma_0$ for $4\leq n\leq 7$.

%\end{enumerate}
%We are going to show that for generic $(A_1,A_2,a_1,a_2)$ there is a solution $(B_1,B_2,b_1,b_2)$ of the equations in \eqref{Barth-matrix} satisfying (i), (ii).

\end{proof}

The lemma implies the rationality of $Y$, and hence the unirationality of $I_n^0$.

\begin{remark}
This method of proof of unirationality does not work for $n\geq 8$, since $\pi|_Y$ is no more dominant. There is a component $\Sigma_0$ of $\Sigma$ dominating $(S^2H_n^\vee)^2\oplus (H_n^\vee)^2$, but 
$\Sigma_0\neq\bar Y$. For $n=8$, both $\Sigma_0$ and $\bar Y$ are of dimension 92, and
$\dim \Sigma_0=n(n+3)+4>\dim \bar Y=8n+\frac12n(n-1)$ for $n\geq 9$. The equality $\dim \Sigma_0=n(n+3)+4$
follows from the fact that if $n\geq 8$, then for generic $(A_i,a_i)$, the space of solutions $(B_i,b_i)$ of the equations in \eqref{Barth-matrix} is 4-dimensional and consists of the quadruples
$$
(B_1,B_2,b_1,b_2)=(\lambda_1 I_n+\lambda_2A_1,\lambda_3 I_n+\lambda_2A_2,
\lambda_4a_1,\lambda_4a_2), \ \ (\lambda_1,\ldots,\lambda_4)\in\CC^4.
$$
\end{remark}

\vspace{0.3cm}
\textbf{Acknowledgements.} A.S.T was supported by the Russian Science Foundation (project 21-41-09011).  D.M. was partially supported by the PRCI SMAGP (ANR-20-CE40-0026-01) and the Labex CEMPI  (ANR-11-LABX-0007-01).


\vspace{2cm}
\begin{thebibliography}{MMMMM}

\bibitem[B1]{B1} {\bf Barth W.},  {\it Moduli of vector bundles on the projective plane},
Inventiones Math. {\bf 42} (1977), 63--91.

\bibitem[B2]{B2} {\bf Barth W.},  {\it Irreducibility of the Space of Mathematical Instanton
Bundles with Rank 2 and $c_2=4$}, Math. Ann. {\bf 258} (1981), 81--106.

\bibitem[BH]{BH} {\bf Barth W., Hulek K.},  {\it Monads and moduli of vector bunbdles}, Manuscripta Math.
{\bf 25} (1978), 323-347.


\bibitem [D]{D} {\bf Dolgachev I. V.}, {\it Rationality of fields of invariants}, In: Algebraic geometry, Bowdoin, 1985, 3--16, Proc. Sympos. Pure Math., 46, Part 2, Amer. Math. Soc., Providence, RI, 1987.


\bibitem [ES]{ES} {\bf Ellingsrud G., Str\o mme S.A.}, {\em Stable rank-2 vector bundles on
$\PP^3$ with $c_1 = 0$
and $c_2 = 3$}, Math. Ann. {\bf 255} (1981), 123--135.



\bibitem [Gr]{Gr} {\bf Grothendieck A.}, {\it Torsion homologique et sections rationnelles},
S\'eminaire Claude Chevalley {\bf 3} (1958), exp. n$^\circ$ 5, pp. 1--29.


\bibitem [H]{H} {\bf Hartshorne R.}, {\em Stable vector bundles of rank 2 on
$\PP^3$}, Math. Ann. {\bf 238} (1978), 229--280.

\bibitem[JV]{JV}
{\bf Jardim M., Verbitsky M.}, Trihyperk\"ahler reduction and instanton 
bundles on $\mathbb{C}\mathbb{P}^3$. Compos. Math. \textbf{150} 
(2014), 1836--1868.


\bibitem [Ka]{Ka} {\bf Katsylo P. I.}, {\em  Rationality of the module variety of mathematical instantons with $c_2=5$}, In: Lie groups, their discrete subgroups, and invariant theory, 105--111, Adv. Soviet Math., 8, Amer. Math. Soc., Providence, RI, 1992.

\bibitem[M2]{M2} 
Daniel R.~Grayson and Michael Stillman, 
Macaulay2, 
A software system for research in algebraic geometry, 
available at 
\url{https://macaulay2.com}.



\bibitem[T]{Tikhomirov2}
{\bf Tikhomirov A.S.}, Moduli of mathematical instanton vector bundles 
with even $c_2$ on projective space. Izvestiya: Mathematics 
\textbf{77}:6 (2013), 1195--1223.


\end{thebibliography}
\end{document}